\documentclass{amsart}
\usepackage{amssymb}
\begin{document}
\newtheorem{theorem}{Theorem}[section]
\newtheorem{lemma}[theorem]{Lemma}
\newtheorem{definition}[theorem]{Definition}
\newtheorem{corollary}[theorem]{Corollary}
\newtheorem{remark}[theorem]{Remark}
\newtheorem{example}[theorem]{Example}
\def\BB{{\mathcal{B}}}
\def\Tr{\operatorname{Tr}}
\def\BOX{\hbox{$\rlap{$\sqcap$}\sqcup$}}
\def\Pspec{\operatorname{Spec}}
\def\vol{\operatorname{Vol}}
\def\II{\mathcal{I}}
\def\binom#1#2{\begin{array}{l}#1\\#2\end{array}}
\makeatletter
  \renewcommand{\theequation}{%
   \thesection.\arabic{equation}}
  \@addtoreset{equation}{section}
\makeatother
\title {The Spectral Geometry of Einstein Manifolds with
Boundary}
\author{JeongHyeong
Park }\thanks{This work was supported by Korea Science and
Engineering Foundation Grant (R05-2003-000-10884-0)}
\begin{address}{JHP Dept. of Computer \& Applied Mathematics, Honam
University, Seobongdong 59, Gwangsanku,  Gwangju, 506-714 South
Korea.\quad Email:jhpark@honam.ac.kr}\end{address}
\begin{abstract}Let $(M,g)$ be a compact Einstein manifold with smooth
boundary. Let $\Delta_{p,\BB}$ be the realization of the $p$ form
valued Laplacian with a suitable boundary condition $\BB$. Let
$\Pspec(\Delta_{p,\BB})$ be the spectrum where each eigenvalue is
repeated according to multiplicity. We show that certain geometric
properties of the boundary may be spectrally characterized in
terms of this data where we fix the Einstein constant.
\end{abstract}
 \maketitle

\section{Introduction}\label{Sect1}

Let $\Delta_p$ be the Laplace-Beltrami operator acting on the
space of smooth $p$ forms over a compact $m$ dimensional
Riemannian manifold $M$ with smooth boundary. If the boundary is
non-empty, then we impose boundary conditions defined by a
suitably chosen operator $\BB$ to define the realization
$\Delta_{p,\BB}$. Let $\nabla$ be the Levi-Civita connection of
$M$ and let $e_m$ be the inward unit normal vector field on the
boundary. Then, for example, Dirichlet and Neumann boundary
conditions are defined by the corresponding Dirichlet and Neumann
boundary operators:
$$\BB_D\phi:=\phi|_{\partial M}
\quad\text{and}\quad\BB_N\phi:=\nabla_{e_m}\phi|_{\partial M}
\quad\text{for}\quad\phi\in C^\infty(\Lambda^pM)\,.$$

In addition to the boundary conditions defined by these operators,
there are also boundary conditions arising from index theory. Near
the boundary, we decompose a differential form
$$\phi=\textstyle\sum_I\phi_Idy^I+\sum_J\psi_Jdx^m\wedge dy^J$$
into tangential and normal components. {\it Absolute boundary
conditions} are then defined by the operator
$$\BB_a\phi:=\{\textstyle\sum_I\partial_m^x\phi_Idy^I\}|_{\partial
M}\oplus\{\textstyle\sum_J\psi_Jdy^J\}|_{\partial M}\,.$$ Dually,
we may use the Hodge $\star$ operator to define the {\it relative
boundary operator} by setting:
$$\BB_r\phi:=\BB_a\star\phi\,.$$

In previous work \cite{P}, we extended a result of Patodi
\cite{Pa70} from the context of closed Riemannian manifolds to the
context of compact Riemannian manifolds with boundary. As we were
interested in determining if the Einstein condition was spectrally
determined, we worked in the context of manifolds of constant
scalar curvature to show:

\begin{theorem}\label{thm1.1} Let $(M_i,g_i)$ be compact Riemannian manifolds
with smooth boundaries and constant scalar curvatures $\tau_i$ for
$i=1,2$. Let $\BB$ define either Dirichlet or Neumann boundary
conditions. Assume
$\Pspec(\Delta_{p,\BB})(M_1)=\Pspec(\Delta_{p,\BB})(M_2)$ for
$p=0,1,2$. Then if $(M_1,g_1)$ is Einstein, then $(M_2,g_2)$ is
Einstein and $\tau_1=\tau_2$.
\end{theorem}

In this paper, instead of studying the geometry of the interior,
we turn our attention to the geometry of the boundary. Motivated
by Theorem \ref{thm1.1}, we shall assume henceforth that the
manifolds under consideration are Einstein and we shall fix the
Einstein constant, or, equivalently, the scalar curvature $\tau$.

We
recall some basic definitions. Let indices $i,j$ range from $1$ to
$m$ and index a local orthonormal frame $\{e_i\}$ for the tangent
bundle of $M$. Near the boundary, we further normalize the frame
and assume that $e_m$ is the inward unit geodesic normal vector
field. Let indices $a,b$ range from $1$ to $m-1$ and index the
induced local orthonormal frame $\{e_a\}$ for the tangent bundle
of the boundary.

We adopt the {\it Einstein convention} and sum
over repeated indices. Let $L$ be the second fundamental form and let $R_{ijkl}$ be the
Riemann curvature tensor. The {\it normalized mean curvature} $\kappa$, the {\it Ricci tensor}
$\rho$, and the {\it scalar curvature} $\tau$ are then given by:
$$\kappa:=L_{aa},\quad\rho_{ij}:=R_{ikkj},\quad\text{and}\quad\tau:=\rho_{ii}\,.$$
Since $M$ is Einstein, $\rho=\lambda g$ where $\lambda$ is the {\it Einstein constant}. This
implies that
$\tau=m\lambda$. Thus fixing the Einstein constant is equivalent to fixing the scalar curvature.

\begin{definition}\label{def1.2}\rm
We say that the boundary of $(M,g)$ is:\begin{enumerate}
\item {\it totally geodesic} if the second fundamental form
vanishes identically. Equivalently, this means that if a geodesic
in $M$ is tangent to the boundary at a single point, then the
geodesic stays in $\partial M$.
\item {\it minimal} if the normalized mean curvature vanishes
identically. Equivalently, this means that the volume of the
boundary is infinitesimally stationary.
\item {\it totally umbillic} if at each point of the boundary, the second fundamental form has only one
eigenvalue; the eigenvalue in question is allowed to vary with the
point of the boundary.
\item {\it strongly totally umbillic} if the the eigenvalue in (3) is independent of the boundary point
chosen.
\end{enumerate}\end{definition}

We can now state the main results of this paper. We first consider
both Dirichlet and Neumann boundary conditions:

\begin{theorem}\label{thm1.3} For $i=1,2$, let $(M_i,g_i)$ be compact Einstein
manifolds with smooth boundaries. Assume that $\tau_1=\tau_2$ and that
\begin{eqnarray*}
&&\Pspec(\Delta_{0,\BB_D})(M_1)=\Pspec(\Delta_{0,\BB_D})(M_2),\quad\text{and}\\
&&\Pspec(\Delta_{0,\BB_N})(M_1)=\Pspec(\Delta_{0,\BB_N})(M_2)
\end{eqnarray*}
where $\BB_D$ and $\BB_N$ define Dirichlet and Neumann boundary conditions, respectively. Then:
\begin{enumerate}
\item If $\partial M_1$ is totally geodesic, then $\partial M_2$ is totally geodesic.
\item If $\partial M_1$ is minimal, then $\partial M_2$ is minimal.
\item If $\partial M_1$ is totally umbillic, then $\partial M_2$ is totally umbillic.
\item If $\partial M_1$ is strongly totally umbillic, then $\partial M_2$ is strongly totally umbillic.
\end{enumerate}
\end{theorem}

In the previous Theorem, we studied two different boundary conditions for the operator
$\Delta_0$. In the next Theorem, we study two different operators, $\Delta_0$ and $\Delta_1$, and
impose either relative or absolute boundary conditions.

\begin{theorem}\label{thm1.4} For $i=1,2$, let $(M_i,g_i)$ be compact Einstein
manifolds with smooth boundaries. Assume that $\tau_1=\tau_2$ and that
\begin{eqnarray*}
&&\Pspec(\Delta_{0,\BB})(M_1)=\Pspec(\Delta_{0,\BB})(M_2),\quad\text{and}\\
&&\Pspec(\Delta_{1,\BB})(M_1)=\Pspec(\Delta_{0,\BB})(M_2)
\end{eqnarray*} where
$\BB$ denotes either relative or absolute boundary conditions.
Then:
\begin{enumerate}
\item If $\partial M_1$ is totally geodesic, then $\partial M_2$ is totally geodesic.
\item If $\partial M_1$ is minimal, then $\partial M_2$ is minimal.
\item If $\partial M_1$ is totally umbillic, then $\partial M_2$ is totally umbillic.
\item If $\partial M_1$ is strongly totally umbillic, then $\partial M_2$ is strongly totally umbillic.
\end{enumerate}
\end{theorem}

Here is a brief outline to the remainder of this paper. In Section
\ref{Sect2}, we review some facts concerning boundary geometry
which we shall need. In Section \ref{Sect3}, we recall some
previous results concerning the heat trace asymptotics. In Section
\ref{Sect4}, we use these results to complete the proof of
Theorems \ref{thm1.3} and \ref{thm1.4}.

\section{The geometry of the boundary}\label{Sect2}

Central to our proof of Theorems \ref{thm1.3} and \ref{thm1.4} is
the following integral characterization of certain geometric
properties. Let $dy$ denote the Riemannian measure on the boundary
and let $dx$ denote the Riemannian measure on the interior. To simplify the notation, let
$$f[M]=\textstyle\int_Mf(x)dx\quad\text{and}\quad f[\partial
M]=\int_{\partial M}f(y)dy$$
 where $f$ is a scalar function.

\begin{theorem}\label{thm2.1} Let $M$ be a compact $m$ dimensional Riemannian
manifold with smooth boundary $\partial M$.
\begin{enumerate} \item $\partial M$ is totally geodesic if and only if
$L_{ab}L_{ab}[\partial M]=0$.
\item $\partial M$ is minimal if and only if $L_{aa}L_{bb}[\partial M]=0$.
\item $\partial M$ is totally umbillic if and only if\\
$\{(m-1)L_{ab}L_{ab}-L_{aa}L_{bb}\}[\partial M]=0$.
\item $\partial M$ is strongly totally umbillic if and only if there exists a
constant $\mu$ so that $\{L_{ab}L_{ab}-2\mu
L_{aa}+\mu^2(m-1)\}[\partial M]=0$.
\end{enumerate}
\end{theorem}

\begin{proof} The first two assertions are immediate. To prove Assertion (3),
we let $\{\kappa_1(y),...,\kappa_{m-1}(y)\}$ be the eigenvalues of
the second fundamental form at a point $y$ of the boundary. Then
the second fundamental form is umbillic at $y$ if and only if
$\kappa_1(y)=...=\kappa_{m-1}(y)$ or equivalently if
$$
0=\textstyle\sum_{i<j}(\kappa_i-\kappa_j)^2\,.
$$
Assertion (3) now follows since we have that
\begin{eqnarray*}
&&\textstyle\sum_i\kappa_i(y)^2=L_{ab}L_{ab}(y),\\
&&\textstyle\sum_{i,j}\kappa_i(y)\kappa_j(y)=L_{aa}L_{bb}(y),\qquad\text{and}\\
&&0\le\textstyle\sum_{i<j}(\kappa_i(y)-\kappa_j(y))^2
=(m-1)L_{ab}(y)L_{ab}(y)-L_{aa}L_{bb}(y)\,.
\end{eqnarray*}
Finally, to prove assertion (4), we note that the second
fundamental form is $\mu$ times the identity at a point $y$ of the
boundary if and only if
$$0=|L-\mu\operatorname{id}|^2=L_{ab}(y)L_{ab}(y)-2\mu L_{aa}(y)+(m-1)\mu^2\,.$$
Since $|L-\mu\operatorname{id}|^2$ is non-negative, Assertion
(4) now holds.
\end{proof}

\section{Heat trace asymptotics}\label{Sect3} To deal with Dirichlet, Neumann,
and absolute boundary conditions in a common framework, it is
useful to introduce the more general notion of {\it mixed boundary
conditions}. Let $\chi$ be a self-adjoint endomorphism of
$\Lambda^p(M)|_{\partial M}$ so that $\chi^2=\operatorname{id}$.
Let $\Pi_\pm$ be orthonormal projection on the $\pm1$ eigenspaces
of $\chi$. Let $S$ be an auxiliary endomorphism of
$\operatorname{range}\Pi_+$. The mixed boundary operator
$\BB_{\chi,S}$ is then defined by
$$\BB_{\chi,S}\phi:=\{\Pi_+(\phi_{;m}+S\phi)\}|_{\partial
M}\oplus\{\Pi_-\phi\}|_{\partial M}\,.$$
\begin{example}\label{exm3.1}
\rm Let $\BB=\BB_{\chi,S}$. \begin{enumerate}
\item If we take $\chi=-\operatorname{id}$, then $\BB$ defines Dirichlet
boundary conditions.
\item If we take $\chi=\operatorname{id}$, then $\BB$ defines Neumann
boundary conditions.
\item Let $\operatorname{ext}(e_i)$ denote left exterior multiplication by the
covector $e_i$ and let $\operatorname{int}(e_i)$ be the dual
operation, left interior multiplication by the covector $e_i$. Let
$\Pi_+$ be projection on $\Lambda(\partial M)$, let $\Pi_-$ be
projection on $\Lambda(\partial M)^\perp$, and let
$$S=-\Pi_+\operatorname{ext}(e_a)\operatorname{int}(e_b)L_{ab}\Pi_+\,.$$
Then $\BB_{\chi,S}$ defines absolute boundary conditions, see, for
example, the discussion in \cite{Gi03}. We note for future
reference that
$$\chi_{;a}=2L_{ab}\{\operatorname{ext}(e_b)\operatorname{int}(e_m)
+\operatorname{ext}(e_m)\operatorname{int}(e_b)\}\,.$$
\end{enumerate}\end{example}

Let $\BB=\BB_{\chi,S}$ and let $e^{-t\Delta_{p,\BB}}$ be the
fundamental solution of the heat equation. The pseudo-differential
calculus established by Seeley \cite{Se,Se69b} shows operator is
of trace class and as $t\downarrow0$ there is a complete
asymptotic expansion with locally computable coefficients in the
form:
$$
\Tr_{L^2}e^{-t\Delta_{p,\BB}}
\sim\textstyle\sum_{n\ge0}t^{(n-m)/2}a_n(\Delta_p,\BB)\,.
$$

Let `;' denote multiple covariant differentiation. The Weitzenb\"och formula permits us to express
$$\Delta_p\omega=-(\omega_{;kk}+E_p\omega)$$
where $E_p$ is a suitably chosen expression in the curvature tensor. For example, we have that
\begin{equation}\label{eqn3.1}
E_0=0\quad\text{and}\quad E_1(e_i)=-\rho_{ij}e_j\,.
\end{equation}
The following result is
a special case of a more general result established by Branson and
Gilkey \cite {BrGi}.

\begin{theorem}\label{thm3.2}
Let $M$ be a compact Riemannian manifold which has a smooth boundary
$\partial M$. Let $\BB=\BB_{\chi,S}$ define mixed boundary
conditions on $\Lambda^p(M)$.
\begin{enumerate}
\smallskip\item
$a_{0}(\Delta_p,\BB)=(4\pi)^{-m/2}\Tr\{\operatorname{id}\}[M]$.
\smallskip\item $a_{1}(\Delta_p,\BB)=(4\pi)^{-(m-1)/2}\frac14
  \Tr\{\chi\}[\partial M]$.
\smallskip\item
$a_{2}(\Delta_p,\BB)=(4\pi)^{-m/2}\frac16\big\{\Tr\{6E_p+\tau\}[M]
+\Tr\{2L_{aa}+12S\}[\partial M]\big\}$.
\smallbreak\item
$a_{3}(\Delta_p,\BB)=(4\pi)^{-(m-1)/2}\frac1{384}
  \Tr\{96
\chi E_p+16 \chi  \tau -8 \chi \rho_{mm}$  \smallbreak\qquad $
+[13 \Pi_+-7 \Pi_-]L_{aa}L_{bb}+[2 \Pi_++10
\Pi_-]L_{ab}L_{ab}+96SL_{aa}$ \smallbreak\qquad $ +192S^{ 2}-12
\chi_{;a} \chi_{;a}\}[\partial M]$.
\end{enumerate}
\end{theorem}

\section{Proof of Theorems \ref{thm1.3} and \ref{thm1.4}}\label{Sect4}

Let $\BB$ denote Dirichlet, Neumann, or absolute boundary
conditions. The heat trace asymptotics $a_n(\Delta_p,\BB)$ are
spectral invariants. Consequently by Theorem \ref{thm3.2},
$$\{\operatorname{vol}(M),\ \operatorname{vol}(\partial M)\}$$
are spectral invariants. We have fixed the Einstein constant and
set the scalar curvature
$\tau=c$. Thus
$$\tau[M],\quad \tau[\partial M],\quad\text{and}\quad \rho_{mm}[\partial M]$$
are spectral invariants as well. The formula for $a_2$ then shows
that
$$L_{aa}[\partial M]$$
is spectrally determined. In light of Theorem \ref{thm2.1}, to
complete the proof of Theorems \ref{thm1.3} and \ref{thm1.4} it
suffices to show
\begin{equation}\label{eqn4.1}
L_{aa}L_{bb}[\partial M]\quad\text{and}\quad L_{ab}L_{ab}[\partial M]
\end{equation}
are spectrally determined by
$\{\Delta_{0,\BB_D},\Delta_{0,\BB_N}\}$, by
$\{\Delta_{0,\BB_a},\Delta_{1,\BB_a}\}$, or by
$\{\Delta_{0,\BB_r},\Delta_{1,\BB_r}\}$.

We shall supress the coefficients of certain invariants in what
follows since they define invariants which are already known to be
spectrally determined; we denote such coefficients by a generic
symbol $\star$. We use the discussion in Example \ref{exm3.1}, the formulae in Equation
\ref{eqn3.1}, and Theorem \ref{thm3.2} to compute:
\begin{eqnarray*}
&&a_3(\Delta_0,\BB_D)=(4\pi)^{(1-m)/2}\textstyle\frac1{384}\{
\star\tau+\star\rho_{mm}-7L_{aa}L_{bb}\\
&&\qquad\qquad\qquad\qquad\qquad\qquad\quad+10L_{ab}L_{ab}\}[\partial M],\\
&&a_3(\Delta_0,\BB_N)=(4\pi)^{(1-m)/2}\textstyle\frac1{384}\{
\star\tau+\star\rho_{mm}+13L_{aa}L_{bb}\\
&&\qquad\qquad\qquad\qquad\qquad\qquad\quad
+2L_{ab}L_{ab}\}[\partial M]\,.
\end{eqnarray*}
Since the coefficient matrix
$$\left(\begin{array}{rr}-7&10\\13&2\end{array}\right)$$
is non-singular, the invariants given in Equation (\ref{eqn4.1})
are spectral invariants as desired; Theorem \ref{thm1.3} now
follows.

To establish Theorem \ref{thm1.4}, we must perform similar
computations for absolute and for relative boundary conditions.

Absolute boundary
conditions are pure Neumann boundary conditions on $0$ forms. By Example \ref{exm3.1},
$$\begin{array}{rr}
\Tr_{\Lambda^1M}\{[13\Pi_+-7\Pi_-]L_{aa}L_{bb}\}=&(13m-20)L_{aa}L_{bb},\\
\Tr_{\Lambda^1M}\{[2\Pi_++10\Pi_-]L_{ab}L_{ab}\}=&(2m+8)L_{ab}L_{ab},\\
\Tr_{\Lambda^1M}\{96SL_{aa}\}=&-96L_{aa}L_{bb},\\
\Tr_{\Lambda^1M}\{192S^2\}=&192L_{ab}L_{ab},\\
\Tr_{\Lambda^1M}\{-12\chi_{:a}\chi_{:a}\}=&-96L_{ab}L_{ab}\,.
\end{array}
$$
It is now an easy matter to use Theorem \ref{thm3.2} to see that
\begin{eqnarray*}
&&a_3(\Delta_0,\BB_a)=(4\pi)^{(1-m)/2}\textstyle\frac1{384}\{
\star\tau+\star\rho_{mm}+13L_{aa}L_{bb}\\
&&\qquad\qquad\qquad\qquad\qquad\qquad\quad
+2L_{ab}L_{ab}\}[\partial M],\\
&&a_3(\Delta_1,\BB_a)=(4\pi)^{(1-m)/2}\textstyle\frac1{384}\{
\star\tau+\star\rho_{mm}+(13m-116)L_{aa}L_{bb}\\
&&\qquad\qquad\qquad\qquad\qquad\qquad\quad+(2m+104)L_{ab}L_{ab}\}[\partial M]\,.
\end{eqnarray*}
The desired result for absolute boundary conditions now follows as
the determinant of the coefficient matrix
$$\left(\begin{array}{rr}13&2\\13m-116&2m+104\end{array}\right)$$
is $1584$ which is different from zero.

Relative boundary conditions are Dirichlet boundary conditions on
$0$ forms. We use duality to see relative boundary conditions on
$1$ forms have the same spectral asymptotics as absolute boundary
conditions on $m-1$ forms. On $m-1$ forms, we compute:
$$\begin{array}{rr}
\Tr_{\Lambda^{m-1}M}\{[13\Pi_+-7\Pi_-]L_{aa}L_{bb}\}=&(-7m+20)L_{aa}L_{bb},\\
\Tr_{\Lambda^{m-1}M}\{[2\Pi_++10\Pi_-]L_{ab}L_{ab}\}=&(10m-8)L_{ab}L_{ab},\\
\Tr_{\Lambda^{m-1}M}\{96SL_{aa}\}=&-96 L_{aa}L_{bb},\\
\Tr_{\Lambda^{m-1}M}\{192S^2\}=&192L_{aa}L_{bb},\\
\Tr_{\Lambda^{m-1}M}\{-12\chi_{:a}\chi_{:a}\}=&-96L_{ab}L_{ab}\,.
\end{array}
$$
It is now an easy matter to use Theorem \ref{thm3.2} to see that
\begin{eqnarray*}
&&a_3(\Delta_0,\BB_r)=(4\pi)^{(1-m)/2}\textstyle\frac1{384}\{
\star\tau+\star\rho_{mm}-7L_{aa}L_{bb}\\
&&\qquad\qquad\qquad\qquad\qquad\qquad\quad+10L_{ab}L_{ab}\}[\partial M],\\
&&a_3(\Delta_1,\BB_r)=(4\pi)^{(1-m)/2}\textstyle\frac1{384}\{
\star\tau+\star\rho_{mm}+(-7m+116)L_{aa}L_{bb}\\
&&\qquad\qquad\qquad\qquad\qquad\qquad\quad+(10m+104)L_{ab}L_{ab}\}[\partial M]\,.
\end{eqnarray*}
The coefficient matrix
$$\left(\begin{array}{rr}-7&10\\-7m+116&10m-104\end{array}\right)$$
has determinant $-432$ which again is different from zero.
\hfill\BOX

\end{document}